\documentclass[11pt,reqno]{amsart}
\usepackage{amsfonts,amssymb,amsmath,color}

\DeclareMathOperator{\bmin}{\mathbf{min}}
\DeclareMathOperator{\bminit}{\boldsymbol{min}}
\DeclareMathOperator{\bmax}{\mathbf{max}}

\DeclareMathOperator{\conv}{\mathrm{conv}}

\newtheorem{theorem}{Theorem}

\newtheorem{corollary}[theorem]{Corollary}

\newtheorem{proposition}[theorem]{Proposition}

\newtheorem{example}[theorem]{Example}

\numberwithin{equation}{section} \numberwithin{theorem}{section}

\begin{document}

\title{Pattern Recognition on Oriented Matroids: Halfspaces, Convex Sets and Tope Committees}

\author{Andrey O. Matveev}
\email{andrey.o.matveev@gmail.com}


\begin{abstract}
The principle of inclusion-exclusion is applied to subsets of maximal covectors contained in halfspaces of a simple oriented matroid and to convex subsets of its ground set for enumerating tope committees.
\end{abstract}

\maketitle

\pagestyle{myheadings}

\markboth{PATTERN RECOGNITION ON ORIENTED MATROIDS}{A.O.~MATVEEV}

\thispagestyle{empty}

\tableofcontents

\section{Introduction}

Let $\mathcal{M}=(E_t,\mathcal{T})$ be a simple oriented matroid (it has no loops, parallel or {\sl antiparallel\/} elements) on the ground set $E_t:=\{1,\ldots,t\}$, with set of topes~$\mathcal{T}$. Throughout we will suppose that $\mathcal{M}$ is not acyclic.

The~family $\mathbf{K}^{\ast}_k(\mathcal{M})$ of {\em tope committees}, of cardinality $k$, $3\leq k\leq|\mathcal{T}|-3$, for the oriented matroid $\mathcal{M}$ is defined as the collection
\begin{equation*}
\mathbf{K}^{\ast}_k(\mathcal{M}):=\left\{\mathcal{K}^{\ast}\subset\mathcal{T}:\ |\mathcal{K}^{\ast}|=k,\
|\mathcal{K}^{\ast}\cap\mathcal{T}_e^+|>
\tfrac{k}{2}\ \ \ \forall e\in E_t\right\}\ ,
\end{equation*}
where $\mathcal{T}_e^+:=\{T\in\mathcal{T}:\ T(e)=+\}$ is the {\em positive halfspace\/} of $\mathcal{M}$ that corresponds to the element $e$, see~\cite{M-Layers,M-Existence,M-Three,M-Reorientations}. The family of {\em tope anti-committees}, of cardinality $k$, for $\mathcal{M}$ is denoted by $\mathbf{A}^{\ast}_k(\mathcal{M})$; by definition, $\mathcal{A}^{\ast}\in\mathbf{A}^{\ast}_k(\mathcal{M})$ iff $-\mathcal{A}^{\ast}\in\mathbf{K}^{\ast}_k(\mathcal{M})$, where~$-\mathcal{A}^{\ast}:=\{-T:\ T\in \mathcal{A}^{\ast}\}$.

Denote by $\binom{\mathcal{T}}{k}$ the family of all $k$-subsets of the tope set $\mathcal{T}$, and consider the families of tope subsets $\mathbf{N}^{\ast}_k(\mathcal{M}):=\binom{\mathcal{T}}{k}-
(\mathbf{K}^{\ast}_k\dot\cup\mathbf{A}^{\ast}_k)$, $3\leq k\leq|\mathcal{T}|-3$, that is the families
\begin{multline*}\mathbf{N}^{\ast}_k(\mathcal{M}):=
\{\mathcal{N}^{\ast}\subset\mathcal{T}:\ |\mathcal{N}^{\ast}|=k,\\
\text{$\mathcal{N}^{\ast}$ neither a committee nor an anti-committee}\}\ ;
\end{multline*}
we have
\begin{equation*}
\#\mathbf{K}^{\ast}_k(\mathcal{M})=
\#\mathbf{A}^{\ast}_k(\mathcal{M})=
\tfrac{1}{2}\bigl(\;\tbinom{|\mathcal{T}|}{k}-\#\mathbf{N}^{\ast}_k(\mathcal{M})\;\bigr)\ ,\ \ \ 3\leq k\leq|\mathcal{T}|-3\ .
\end{equation*}

For an element $e\in E_t$, we let $\mathcal{T}_e^-:=\{T\in\mathcal{T}:\ T(e)=-\}$ denote the {\em negative halfspace\/} of $\mathcal{M}$ that corresponds to the element $e$.
The family of all subsets, of cardinality $j$, of the positive halfspace~$\mathcal{T}_e^+$ is denoted by~$\tbinom{\mathcal{T}_e^+}{j}$ and, similarly,~$\tbinom{\mathcal{T}_e^-}{i}$ denotes the family of $i$-subsets of the negative halfspace~$\mathcal{T}_e^-$. The family of $(i+j)$-sets $\tbinom{\mathcal{T}_e^-}{i}\boxplus\tbinom{\mathcal{T}_e^+}{j}$ is defined as the family~$\bigl\{A\dot\cup B:\ A$ $\in\tbinom{\mathcal{T}_e^-}{i},\ B\in\tbinom{\mathcal{T}_e^+}{j}\bigr\}$.

On the one hand, $\mathbf{K}^{\ast}_k(\mathcal{M})=\bigcap_{e\in E_t}\bigcup_{\lceil(k+1)/2\rceil\leq j\leq k}\left(\tbinom{\mathcal{T}_e^-}{k-j}\boxplus\tbinom{\mathcal{T}_e^+}{j}\right)$, $3\leq k$ $\leq |\mathcal{T}|-3$.
On the other hand, a $k$-subset $\mathcal{K}^{\ast}\subset\mathcal{T}$ is a committee for $\mathcal{M}$ iff
\begin{itemize}
\item the set $\mathcal{K}^{\ast}$ contains no set from the family $\bigcup_{e\in E_t}\binom{\mathcal{T}^-_{e}}{\lfloor(k+1)/2\rfloor}$;
\item the set $\mathcal{K}^{\ast}$ contains at least one set from each family $\binom{\mathcal{T}^+_{e}}{\lceil(k+1)/2\rceil}$, $e\in E_t$ --- in other words,
    the collection $\binom{\mathcal{K}^{\ast}}{\lceil(k+1)/2\rceil}$ is a {\em blocking family} for the family
$\{\binom{\mathcal{T}_1^+}{\lceil(k+1)/2\rceil},\ldots,\binom{\mathcal{T}_t^+}{\lceil(k+1)/2\rceil}\}$, that is, $\#(\binom{\mathcal{K}^{\ast}}{\lceil(k+1)/2\rceil}$ $\cap\binom{\mathcal{T}_e^+}{\lceil(k+1)/2\rceil})$ $>0$, $e\in E_t$.
\end{itemize}
As a consequence, the collection $\mathbf{K}^{\ast}_k(\mathcal{M})$ is the family of all {\em blocking $k$-sets\/} of topes for the family $\bigcup_{e\in E_t}\binom{\mathcal{T}^+_e}{\lfloor(|\mathcal{T}|-k+1)/2\rfloor}$, and a committee $\mathcal{K}^{\ast}\in\mathbf{K}^{\ast}_k(\mathcal{M})$ is {\em minimal\/} if any its proper $i$-subset
$\mathcal{I}^{\ast}\subset\mathcal{K}^{\ast}$ is not a blocking set for the family $\bigcup_{e\in E_t}\binom{\mathcal{T}^+_e}{\lfloor(|\mathcal{T}|-i+1)/2\rfloor}$.

Based on these remarks, we calculate in Sections~\ref{s:4} and~\ref{s:2} the numbers~$\#\mathbf{K}^{\ast}_k(\mathcal{M})$ of general committees of cardinality $k$, in several possible ways, by applying the {\em principle of inclusion-exclusion\/} \cite{A,St1} to subsets of maximal covectors  contained in halfspaces of the oriented matroid $\mathcal{M}$; in Section~\ref{s:2}, these calculations involve the {\em convex subsets\/} of the ground set of $\mathcal{M}$. In~Section~\ref{s:3} we find the numbers~$\#\overset{\circ}{\mathbf{K}}{}^{\ast}_k(\mathcal{M})$ of tope committees, of cardinality $k$, which contain no pairs of opposites. Sections~\ref{s:1}, \ref{s:5} and~\ref{s:6} list auxiliary results.

See~\cite{ERW} and references therein on acyclic, convex and free sets of oriented matroids.

One can associate to the oriented matroid $\mathcal{M}$ various ``{\sl $\kappa^{\ast}$\!-vectors}'' (and their flag generalizations) whose components are the numbers of its tope committees of the corresponding cardinality, for example:
\begin{itemize}
\item the vector $\pmb{\kappa}^{\ast}(\mathcal{M}):=\bigl(\kappa^{\ast}_1(\mathcal{M}),\ldots,
    \kappa^{\ast}_{|\mathcal{T}|/2}(\mathcal{M})\bigr)\in\mathbb{N}^{|\mathcal{T}|/2}$, where
$\kappa^{\ast}_k(\mathcal{M})$ $:=$\hfill $\#\mathbf{K}^{\ast}_k(\mathcal{M})$\hfill  ---\hfill  recall\hfill  that\hfill  $\#\mathbf{K}^{\ast}_k(\mathcal{M})$\hfill  $=$\hfill  $\#\mathbf{K}^{\ast}_{|\mathcal{T}|-k}(\mathcal{M})$, \hfill $1\leq k$\\ $\leq
|\mathcal{T}|-1$;

\item the vector
$\overset{\circ}{\pmb{\kappa}}{}^{\ast}(\mathcal{M}):=\bigl(\overset{\circ}{\kappa}{}^{\ast}_1(\mathcal{M}),\ldots,
    \overset{\circ}{\kappa}{}^{\ast}_{|\mathcal{T}|/2}(\mathcal{M})\bigr)\in\mathbb{N}^{|\mathcal{T}|/2}$, where
$\overset{\circ}{\kappa}{}^{\ast}_k(\mathcal{M})$ $:=$\hfill  $\#\overset{\circ}{\mathbf{K}}{}^{\ast}_k(\mathcal{M})$\hfill  ---\hfill  note\hfill  that \hfill $\#\overset{\circ}{\mathbf{K}}{}^{\ast}_k(\mathcal{M})=0$\hfill  whenever\hfill  $\frac{1}{2}|\mathcal{T}|< k$\\ $\leq
|\mathcal{T}|-1$;

\item the vector $\pmb{\kappa}^{\ast}_{\min}(\mathcal{M}):=\bigl(\kappa^{\ast}_{\min\; 1}(\mathcal{M}),\ldots,
    \kappa^{\ast}_{\min\; |\mathcal{T}|/2}(\mathcal{M})\bigr)\in\mathbb{N}^{|\mathcal{T}|/2}$, where
$\kappa^{\ast}_{\min\; k}(\mathcal{M}):=\#\{\mathcal{K}^{\ast}\in\mathbf{K}^{\ast}_k(\mathcal{M})$ $:\ \text{$\mathcal{K}^{\ast}$ minimal}\}$;

\item the vector $\pmb{\kappa}^{\ast}_{\bmax^+}(\mathcal{M}):=\bigl(\kappa^{\ast}_{\bmax^+\; 1}(\mathcal{M}),\ldots,
    \kappa^{\ast}_{\bmax^+\; |\mathcal{T}|/2}(\mathcal{M})\bigr)\in\mathbb{N}^{|\mathcal{T}|/2}$, where $\kappa^{\ast}_{\bmax^+\; k}(\mathcal{M}):=\#\{\mathcal{K}^{\ast}\in\mathbf{K}^{\ast}_k(\mathcal{M}):\ \mathcal{K}^{\ast}
\subseteq\bmax^+(\mathcal{T})\}$, the number of tope committees, of cardinality $k$, composed of topes from the set $\bmax^+(\mathcal{T})\subset\mathcal{T}$ of all topes of $\mathcal{M}$ with inclusion-maximal positive parts,
\end{itemize}
 and so on.

\begin{example} Let $\mathcal{M}:=(E_6,\mathcal{T})$ be the
simple oriented matroid represented by its third positive halfspace
{\tiny
\begin{equation*}
\mathcal{T}^+_3:=
\begin{matrix}
\{&-&-&+&+&+&+\\ &-&-&+&-&+&+\\ &+&-&+&-&+&+\\ &+&-&+&-&+&-&\\
&-&-&+&-&+&-\\ &-&-&+&+&+&-\\ &-&-&+&+&-&+\\ &-&+&+&+&-&+&\\
&-&+&+&+&+&-\\ &-&+&+&-&+&-\\ &+&+&+&-&+&-\\ &+&+&+&-&-&-&\\
&-&+&+&-&-&-\\ &-&+&+&+&-&-&\}\ ;
\end{matrix}
\end{equation*}
}

\noindent a realization of its reorientation ${}_{-\{1,2\}}\mathcal{M}$, by a hyperplane arrangement in~$\mathbb{R}^3$, is shown in~\cite[Figure~3.1]{M-Existence}.

The oriented matroid $\mathcal{M}$ has
$28$ maximal covectors and $238012$ tope committees ---
\begin{gather*}
\pmb{\kappa}^{\ast}(\mathcal{M})=(\underset{\overset{\uparrow}{1}}{0},
\overset{\underset{\downarrow}{2}}{0},
\underset{\overset{\uparrow}{3}}{3},
\overset{\underset{\downarrow}{4}}{0},\underset{\overset{\uparrow}{5}}{144},\overset{\underset{\downarrow}{6}}{1},
\underset{\overset{\uparrow}{7}}{1942},\overset{\underset{\downarrow}{8}}{22},
\underset{\overset{\uparrow}{9}}{11872},\overset{\underset{\downarrow}{10}}{136},
\underset{\overset{\uparrow}{11}}{37775},\overset{\underset{\downarrow}{12}}{386},
\underset{\overset{\uparrow}{13}}{66454},\overset{\underset{\downarrow}{14}}{542})
\intertext{--- among which $4496$ committees are free of opposites:}
\overset{\circ}{\pmb{\kappa}}{}^{\ast}(\mathcal{M})=(\underset{\overset{\uparrow}{1}}{0},
\overset{\underset{\downarrow}{2}}{0},
\underset{\overset{\uparrow}{3}}{3},
\overset{\underset{\downarrow}{4}}{0},\underset{\overset{\uparrow}{5}}{111},\overset{\underset{\downarrow}{6}}{1},
\underset{\overset{\uparrow}{7}}{778},\overset{\underset{\downarrow}{8}}{14},
\underset{\overset{\uparrow}{9}}{1935},\overset{\underset{\downarrow}{10}}{24},
\underset{\overset{\uparrow}{11}}{1448},\overset{\underset{\downarrow}{12}}{24},
\underset{\overset{\uparrow}{13}}{158},\overset{\underset{\downarrow}{14}}{0})\ .
\end{gather*}
\end{example}

\section{Relative Blocking in Boolean Lattices}
\label{s:1}

Let $\varLambda$ be a nontrivial antichain in the Boolean lattice $\mathbb{B}(n)$ of rank $n$, and~$\varLambda^{\bot}$ the set of lattice complements of the elements of $\varLambda$ in $\mathbb{B}(n)$;
$\rho(\cdot)$ denotes the rank function, $\mathbb{B}(\mathcal{T})^{(i)}:=\{b\in\mathbb{B}(n):\ \rho(b)=i\}$ denotes the~$i$th layer of $\mathbb{B}(n)$,
and $\mathfrak{I}(C)$ stands for the order ideal of the lattice $\mathbb{B}(n)$ generated by its antichain $C$.

For a rational number $r$, $0\leq r<1$, and for a positive integer number $k$, consider the subset
\begin{equation}
\label{eq:12}
\mathbf{I}_{r,k}\bigl(\mathbb{B}(n),
\varLambda\bigr):=\bigl\{b\in\mathbb{B}(n):\ \rho(b)=k,\
\rho(b\wedge\lambda)>r\cdot k\ \ \ \forall
\lambda\in\varLambda\bigr\}\subset\mathbb{B}(n)^{(k)}
\end{equation}
that consists of the {\em relatively $r$-blocking elements}, of rank $k$, for the antichain~$\varLambda$.

Set $\nu(r\cdot k):=\lfloor r\cdot k\rfloor+1$ and consider an antichain $\varLambda\subset\mathbb{B}(n)$ such that~$\rho(\lambda)\geq \nu(r\cdot k)$ and $n-\rho(\lambda)\geq k-\nu(r\cdot k)+1$,
for each element $\lambda\in\varLambda$, that is,
\begin{equation}
\label{eq:22}
\lfloor r\cdot k\rfloor+1\leq \min_{\lambda\in\varLambda}\rho(\lambda)\ \ \ \text{and}\ \ \
\max_{\lambda\in\varLambda}\rho(\lambda)\leq n+\lfloor r\cdot k\rfloor-k\ .
\end{equation}

If\hfill the\hfill antichain\hfill $\varLambda$\hfill satisfies\hfill constraints\hfill (\ref{eq:22})\hfill  then\hfill for\hfill an\hfill element\hfill $b'$\\ $\in\mathbb{B}(n)^{(k)}$ we have $b'\not\in\mathbf{I}_{r,k}\bigl(\mathbb{B}(n),
\varLambda\bigr)$ iff $b'>d'$ for at least one element $d'$ of rank~$k-\nu(r\cdot k)+1=k-\lfloor r\cdot k\rfloor$ such that $d'\in\mathfrak{I}(\varLambda^{\bot})$; therefore, on the one hand,
\begin{multline}
\label{eq:13}
\bigl|\mathbf{I}_{r,k}\bigl(\mathbb{B}(n),\varLambda\bigr)\bigr|=\binom{n}{k}\\ +
\sum_{D'\subseteq \mathbb{B}(n)^{(k-\lfloor r\cdot k\rfloor)}\cap\mathfrak{I}(\varLambda^{\bot}):\
|D'|>0}(-1)^{|D'|}\cdot\binom{n-\rho\left(\bigvee_{d'\in
D'}d'\right)}{n-k}\ .
\end{multline}
On the other hand, for an element $b\in\mathbb{B}(n)$ the inclusion $b\in\mathbf{I}_{r,k}\bigl(\mathbb{B}(n),\varLambda\bigr)$ holds\hfill iff\hfill for\hfill each\hfill element\hfill $\lambda\in\varLambda$ we\hfill have\hfill $\rho(b\wedge\theta_\lambda)>0$,\hfill for\hfill any\hfill element\hfill $\theta_\lambda$\\ $\in\mathbb{B}(n)^{(\rho(\lambda)-\nu(r\cdot k)+1)}\cap\mathfrak{I}(\lambda)$, that is,
\begin{equation*}
b\in\mathbf{I}_{r,k}\bigl(\mathbb{B}(n),\varLambda\bigr)\ \ \ \Longleftrightarrow\ \ \
\rho(b\wedge\theta_\lambda)>0\ \ \ \forall\theta_\lambda\in\mathbb{B}(n)^{(\rho(\lambda)-\lfloor r\cdot k\rfloor)}\cap\mathfrak{I}(\lambda)\ \ \forall\lambda\in\varLambda\ ,
\end{equation*}
and we have
\begin{multline}
\label{eq:15}
\bigl|\mathbf{I}_{r,k}\bigl(\mathbb{B}(n),\varLambda\bigr)\bigr|=\binom{n}{k}\\+
\sum_{D\subseteq\bmin\bigcup_{\lambda\in\varLambda}(
\mathbb{B}(n)^{(\rho(\lambda)-\lfloor r\cdot k\rfloor)}\cap\mathfrak{I}(\lambda)):\ |D|>0}
(-1)^{|D|}\cdot\binom{n-\rho(\bigvee_{d\in D}d)}{k}
\end{multline}
(where $\bmin\cdot$ stands for the set of minimal elements of a subposet) or, via {\em Vandermonde's convolution},
\begin{multline}
\label{eq:16}
\bigl|\mathbf{I}_{r,k}\bigl(\mathbb{B}(n),\varLambda\bigr)\bigr|=
-\sum_{D\subseteq\bmin\bigcup_{\lambda\in\varLambda}
(\mathbb{B}(n)^{(\rho(\lambda)-\lfloor r\cdot k\rfloor)}\cap\mathfrak{I}(\lambda)):\ |D|>0}
(-1)^{|D|}\\ \cdot\sum_{1\leq h\leq k}\binom{\rho(\bigvee_{d\in D}d)}{h}\!\!\binom{n-\rho(\bigvee_{d\in D}d)}{k-h}\ .
\end{multline}

One\hfill  more\hfill  inclusion-exclusion type\hfill  formula\hfill  for\hfill  the\hfill  cardinality\hfill  of\hfill  the\\ set~$\mathbf{I}_{r,k}\bigl(\mathbb{B}(n),\varLambda\bigr)$, for an antichain $\varLambda$ such that $\rho(\lambda)\geq\nu(r\cdot k)$, for all $\lambda\in\varLambda$, is given in~\cite[(5.4)]{M-Relative}: if
\begin{equation}
\label{eq:23}
\lfloor r\cdot k\rfloor+1\leq \min_{\lambda\in\varLambda}\rho(\lambda)
\end{equation}
then
\begin{multline}
\label{eq:17}
\bigl|\mathbf{I}_{r,k}\bigl(\mathbb{B}(n),\varLambda\bigr)\bigr|
=\sum_{D\subseteq \mathbb{B}(n)^{(\lfloor r\cdot k\rfloor+1)}\cap\mathfrak{I}(\varLambda):\ |D|>0}
(-1)^{|D|}\;\\ \cdot\left(\sum_{C\subseteq\varLambda:\ D\subseteq\mathfrak{I}(C)}(-1)^{|C|}\right)\!\!\!
\binom{n-\rho(\vee_{d\in D}d)}{n-k}\ .
\end{multline}

We now refine formulas~(\ref{eq:13}), (\ref{eq:15}) and~(\ref{eq:16}) with the help of the {\em M\"{o}bius function}~\cite{A,St1}, see below expressions~(\ref{eq:28}), (\ref{eq:29}) and~(\ref{eq:30}), respectively. Let~$X$ be a nontrivial antichain in the Boolean lattice $\mathbb{B}(n)$. Denote by~$\mathcal{E}(\mathbb{B}(n),$ $X)$ the sub-join-semilattice of $\mathbb{B}(n)$ generated by the set $X$ and augmented by a new least element $\hat{0}$; the greatest element $\hat{1}$ of the lattice~$\mathcal{E}\bigl(\mathbb{B}(n),X\bigr)$ is the join $\bigvee_{x\in X}x$ in $\mathbb{B}(n)$. The M\"{o}bius function of the lattice $\mathcal{E}\bigl(\mathbb{B}(n),X\bigr)$ is denoted by $\mu_{\mathcal{E}}(\cdot,\cdot)$.

Let $\varLambda$ be a nontrivial antichain in the Boolean lattice $\mathbb{B}(n)$ that complies with constraints~(\ref{eq:22}). We have
\begin{itemize}
\item
\begin{multline}
\label{eq:28}
\bigl|\mathbf{I}_{r,k}\bigl(\mathbb{B}(n),\varLambda\bigr)\bigr|=\binom{n}{k}\\ +
\sum_{z\in\mathcal{E}(\mathbb{B}(n)^{(k-\lfloor r\cdot k\rfloor)}\cap\mathfrak{I}(\varLambda^{\bot})):\ z>\hat{0}}\mu_{\mathcal{E}}(\hat{0},z)\cdot\binom{n-\rho(z)}{n-k}\ .
\end{multline}
\item
\begin{multline}
\label{eq:29}
\bigl|\mathbf{I}_{r,k}\bigl(\mathbb{B}(n),\varLambda\bigr)\bigr|=\binom{n}{k}\\+
\sum_{z\in\mathcal{E}(\bmin\bigcup_{\lambda\in\varLambda}(
\mathbb{B}(n)^{(\rho(\lambda)-\lfloor r\cdot k\rfloor)}\cap\mathfrak{I}(\lambda))):\ z>\hat{0}}\mu_{\mathcal{E}}(\hat{0},z)\cdot\binom{n-\rho(z)}{k}\ .
\end{multline}
\item
\begin{multline}
\label{eq:30}
\bigl|\mathbf{I}_{r,k}\bigl(\mathbb{B}(n),\varLambda\bigr)\bigr|=
-\sum_{z\in\mathcal{E}(\bmin\bigcup_{\lambda\in\varLambda}(
\mathbb{B}(n)^{(\rho(\lambda)-\lfloor r\cdot k\rfloor)}\cap\mathfrak{I}(\lambda))):\ z>\hat{0}}\mu_{\mathcal{E}}(\hat{0},z)\\ \cdot\sum_{1\leq h\leq k}\binom{\rho(z)}{h}\!\!\binom{n-\rho(z)}{k-h}\ .
\end{multline}
\end{itemize}

A companion formula to (\ref{eq:17}) is given in~\cite[(5.6)]{M-Relative}: let $\varLambda\subset\mathbb{B}(n)$ be an antichain that obeys constraint~(\ref{eq:23}), and let $\mathcal{C}_{r,k}\bigl(\mathbb{B}(n),\varLambda\bigr)$ be the join-semilattice of all sets
from the family $\{\mathbb{B}(n)^{(\lfloor r\cdot k\rfloor+1)}\cap\mathfrak{I}(C):\ C\subseteq \varLambda,\ |C|>0\}$
ordered by inclusion and augmented by a new least element $\hat{0}$; the greatest element $\hat{1}$ of the lattice $\mathcal{C}_{r,k}\bigl(\mathbb{B}(n),\varLambda\bigr)$ is the set $\mathbb{B}(n)^{(\lfloor r\cdot k\rfloor+1)}\cap\mathfrak{I}(\varLambda)$. We denote the M\"{o}bius function of $\mathcal{C}_{r,k}\bigl(\mathbb{B}(n),\varLambda\bigr)$ by $\mu_{\mathcal{C}}(\cdot,\cdot)$. We have
\begin{multline}
\label{eq:19}
\bigl|\mathbf{I}_{r,k}\bigl(\mathbb{B}(n),\varLambda\bigr)\bigr|=\sum_{X\in\mathcal{C}_{r,k}(\mathbb{B}(n),\varLambda):\ X>\hat{0}}
\mu_{\mathcal{C}}(\hat{0},X)\\
\cdot\sum_{z\in\mathcal{E}(\mathbb{B}(n),X):\ z>\hat{0}}\mu_{\mathcal{E}}(\hat{0},z)\cdot\binom{n-\rho(z)}{n-k}\ .
\end{multline}

\section{Halfspaces and Tope Committees}

\label{s:4}

Let $\mathbb{B}(\mathcal{T})$ be the Boolean lattice of all subsets of the tope set $\mathcal{T}$, and~$\Upsilon$ $:=\{\upsilon_1,\ldots,\upsilon_t\}\subset\mathbb{B}(\mathcal{T})^{(|\mathcal{T}|/2)}$ its antichain whose element $\upsilon_e$ represents in~$\mathbb{B}(\mathcal{T})$ the $e$th positive halfspace $\mathcal{T}^+_e$ of the oriented matroid $\mathcal{M}$.
The family $\mathbf{K}^{\ast}_k(\mathcal{M})$ of tope committees, of cardinality $k$, $3\leq k\leq|\mathcal{T}|-3$, for $\mathcal{M}$ is represented in the lattice $\mathbb{B}(\mathcal{T})$ by the antichain
\begin{equation*}
\mathbf{I}_{\frac{1}{2},k}\bigl(\mathbb{B}(\mathcal{T}),
\Upsilon\bigr):=\bigl\{b\in\mathbb{B}(\mathcal{T}):\ \rho(b)=k,\
\rho(b\wedge\upsilon_e)>\tfrac{k}{2}\ \ \ \forall
e\in E_t \bigr\}\subset\mathbb{B}(\mathcal{T})^{(k)}\ ;
\end{equation*}
thanks to axiomatic symmetry $\mathcal{T}=-\mathcal{T}$,
see~\cite[\S4.1.1,~(L1)]{BLSWZ},
the cardinality of this set is
\begin{multline}
\label{eq:6}
\bigl|\mathbf{I}_{\frac{1}{2},k}\bigl(\mathbb{B}(\mathcal{T}),\Upsilon\bigr)\bigr|=\binom{|\mathcal{T}|}{k}\\ +
\sum_{D\subseteq \mathbb{B}(\mathcal{T})^{(\lfloor(k+1)/2\rfloor)}\cap\mathfrak{I}(\Upsilon):\
|D|>0}(-1)^{|D|}\cdot\binom{|\mathcal{T}|-\rho\left(\bigvee_{d\in
D}d\right)}{|\mathcal{T}|-k}\ ,
\end{multline}
by~(\ref{eq:13}). Note that for an integer $j$, $1\leq j\leq|\mathcal{T}|/2$, we have
\begin{equation*}
\bigl|\mathbb{B}(\mathcal{T})^{(j)}\cap\mathfrak{I}(\Upsilon)\bigr|=-
\sum_{\substack{A\in
L_{\conv}(\mathcal{M})-\{\hat{0}\}
:\\ \text{$A$ free}}}(-1)^{|A|}\cdot\binom{|\mathcal{T}^+_A|}{j}\ ,
\end{equation*}
where $L_{\conv}(\mathcal{M})$ denotes the meet-semilattice of {\em convex subsets\/} of the ground set $E_t$, and $\mathcal{T}^+_A:=\bigcap_{a\in A}\mathcal{T}^+_a$; $\hat{0}$ denotes the least element of~$L_{\conv}(\mathcal{M})$. Recall that from the algebraic combinatorial point of view~\cite{BI}, the set $\mathbb{B}(\mathcal{T})^{(j)}\cap\mathfrak{I}(\Upsilon)$ is a {\em subset\/} in the {\em Johnson association scheme\/} $\mathbf{J}(|\mathcal{T}|,j):=(\pmb{X},\boldsymbol{\mathcal{R}})$
on the set~$\pmb{X}:=\mathbb{B}(\mathcal{T})^{(j)}$,
with the partition
$\boldsymbol{\mathcal{R}}:=(\pmb{R}_0,\pmb{R}_1,\ldots,$
$\pmb{R}_j)$ of $\pmb{X}\times\pmb{X}$,
defined by $\pmb{R}_i:=\bigl\{(x,y):\
j-\rho(x\wedge y)=i\bigr\}$,
for all $0\leq i\leq j$.

Reformulate observation~(\ref{eq:6}) in the following way:
\begin{multline}
\label{eq:5}
\#\mathbf{K}^{\ast}_k(\mathcal{M})=\#\mathbf{K}^{\ast}_{|\mathcal{T}|-k}(\mathcal{M})\\=
\binom{|\mathcal{T}|}{|\mathcal{T}|-\ell}+
\sum_{\substack{\mathcal{G}\subseteq\bigcup_{e\in E_t}\binom{\mathcal{T}^+_e}{\lfloor(\ell+1)/2\rfloor}:\\
1\leq\#\mathcal{G}\leq\binom{\ell}{\lfloor(\ell+1)/2\rfloor},\\
|\bigcup_{G\in\mathcal{G}}G|\leq\ell}}
(-1)^{\#\mathcal{G}}\cdot
\binom{|\mathcal{T}|-|\bigcup_{G\in\mathcal{G}}G|}{|\mathcal{T}|-\ell}\ ,
\end{multline}
where $\ell\in\{k,|\mathcal{T}|-k\}$; this formula counts
the number of all blocking $k$-sets of topes for the family $\bigcup_{e\in E_t}\binom{\mathcal{T}^+_e}{\lfloor(|\mathcal{T}|-k+1)/2\rfloor}$, cf.~(\ref{eq:15}), and it counts
the number of all blocking $(|\mathcal{T}|-k)$-sets of topes for the family $\bigcup_{e\in E_t}\binom{\mathcal{T}^+_e}{\lfloor(k+1)/2\rfloor}$. We can also rewrite~(\ref{eq:6}) by means of Vandermonde's convolution in the form:
\begin{multline*}
\bigl|\mathbf{I}_{\frac{1}{2},k}\bigl(\mathbb{B}(\mathcal{T}),\Upsilon\bigr)\bigr|=
-
\sum_{D\subseteq\mathbb{B}(\mathcal{T})^{(\lfloor(|\mathcal{T}|-k+1)/2\rfloor)}\cap\mathfrak{I}(\Upsilon):\ |D|>0}
(-1)^{|D|}\\ \cdot\sum_{1\leq h\leq k}\binom{\rho(\bigvee_{d\in D}d)}{h}\!\!\binom{|\mathcal{T}|-\rho(\bigvee_{d\in D}d)}{k-h}\ ,
\end{multline*}
cf.~(\ref{eq:16}), that is,
\begin{multline}
\label{eq:24}
\#\mathbf{K}^{\ast}_k(\mathcal{M})=\#\mathbf{K}^{\ast}_{|\mathcal{T}|-k}(\mathcal{M})=
-\sum_{\mathcal{G}\subseteq\bigcup_{e\in E_t}\binom{\mathcal{T}^+_e}{\lfloor(|\mathcal{T}|-\ell+1)/2\rfloor}:\
\#\mathcal{G}>0}
(-1)^{\#\mathcal{G}}\\ \cdot\sum_{\max\{1,\ell-|\mathcal{T}|+|\bigcup_{G\in\mathcal{G}}G|\}\leq h\leq \min\{\ell,|\bigcup_{G\in\mathcal{G}}G|\}}\binom{|\bigcup_{G\in\mathcal{G}}G|}{h}
\!\!\binom{|\mathcal{T}|-|\bigcup_{G\in\mathcal{G}}G|}{\ell-h}\ ,
\end{multline}
where $\ell\in\{k,|\mathcal{T}|-k\}$.

If $\mathcal{G}$ is a family of tope subsets then we denote by $\boldsymbol{\mathcal{E}}(\mathcal{G})$ the join-semilattice $\{\bigcup_{F\in\mathcal{F}}F:\ \mathcal{F}\subseteq\mathcal{G},\ \#\mathcal{F}>0\}$ composed of the unions of the sets from the family $\mathcal{G}$ ordered by inclusion and augmented by a new least element $\hat{0}$; the greatest element $\hat{1}$ of the lattice $\boldsymbol{\mathcal{E}}(\mathcal{G})$ is the set $\bigcup_{G\in\mathcal{G}}G$. The M\"{o}bius function of the poset $\boldsymbol{\mathcal{E}}(\mathcal{G})$ is denoted by $\mu_{\boldsymbol{\mathcal{E}}}(\cdot,\cdot)$.

Expressions~(\ref{eq:25}) and (\ref{eq:26}) below refine formulas~(\ref{eq:5}) and~(\ref{eq:24}), respectively.

\begin{proposition}
The number $\#\mathbf{K}^{\ast}_k(\mathcal{M})$ of tope committees, of cardinality~$k$, $3\leq k\leq
|\mathcal{T}|-3$, for the oriented matroid $\mathcal{M}:=(E_t,\mathcal{T})$, is:

{\rm(i)}

\nopagebreak[4]

\begin{multline}
\label{eq:25}
\#\mathbf{K}^{\ast}_k(\mathcal{M})=\#\mathbf{K}^{\ast}_{|\mathcal{T}|-k}(\mathcal{M})\\
=\binom{|\mathcal{T}|}{|\mathcal{T}|-\ell} +
\sum_{G\in\boldsymbol{\mathcal{E}}(\bigcup_{e\in E_t}\binom{\mathcal{T}^+_e}{\lfloor(\ell+1)/2\rfloor}):\ 0<|G|\leq\ell}\mu_{\boldsymbol{\mathcal{E}}}(\hat{0},G)\cdot\binom{|\mathcal{T}|-|G|}{|\mathcal{T}|-\ell}\ ,
\end{multline}
where $\ell\in\{k,|\mathcal{T}|-k\}$.

{\rm(ii)}

\begin{multline}
\label{eq:26}
\#\mathbf{K}^{\ast}_k(\mathcal{M})=\#\mathbf{K}^{\ast}_{|\mathcal{T}|-k}(\mathcal{M})=
-\sum_{G\in\boldsymbol{\mathcal{E}}(\bigcup_{e\in E_t}\binom{\mathcal{T}^+_e}{\lfloor(|\mathcal{T}|-\ell+1)/2\rfloor}):\ |G|>0}\mu_{\boldsymbol{\mathcal{E}}}(\hat{0},G)\\ \cdot\sum_{\max\{1,\ell-|\mathcal{T}|+|G|\}\leq h\leq \min\{\ell,|G|\}}\binom{|G|}{h}
\!\!\binom{|\mathcal{T}|-|G|}{\ell-h}\ ,
\end{multline}
where $\ell\in\{k,|\mathcal{T}|-k\}$.
\end{proposition}

Let\hfill $\mathcal{C}_{\frac{1}{2},k}\bigl(\mathbb{B}(\mathcal{T}),\Upsilon\bigr)$\hfill be\hfill
the\hfill join-semilattice\hfill of\hfill all\hfill sets\hfill
from\hfill the\hfill family\\ $\{\mathbb{B}(\mathcal{T})^{(\lceil(k+1)/2\rceil)}\cap\mathfrak{I}(C):\ C\subseteq \Upsilon,\ |C|>0\}$
ordered by inclusion and augmented by a new least element $\hat{0}$. The greatest element $\hat{1}$ of the lattice~$\mathcal{C}_{\frac{1}{2},k}\bigl(\mathbb{B}(\mathcal{T}),\Upsilon\bigr)$ is the set $\mathbb{B}(\mathcal{T})^{(\lceil(k+1)/2\rceil)}\cap\mathfrak{I}(\Upsilon)$. Similarly,
for an element $X\in\mathcal{C}_{\frac{1}{2},k}\bigl(\mathbb{B}(\mathcal{T}),\Upsilon\bigr)$ we denote by $\mathcal{E}\bigl(\mathbb{B}(\mathcal{T}),X\bigr)$ the sub-join-semilattice of $\mathbb{B}(\mathcal{T})$ generated by the set $X\subset\mathbb{B}(\mathcal{T})$ and augmented by a new least element $\hat{0}$.
The M\"{o}bius functions of the posets $\mathcal{C}_{\frac{1}{2},k}\bigl(\mathbb{B}(\mathcal{T}),\Upsilon\bigr)$ and $\mathcal{E}\bigl(\mathbb{B}(\mathcal{T}),X\bigr)$ are denoted by $\mu_{\mathcal{C}}(\cdot,\cdot)$ and $\mu_{\mathcal{E}}(\cdot,\cdot)$, respectively.

Using~(\ref{eq:19}), we obtain the expression
\begin{multline}
\label{eq:27}
\bigl|\mathbf{I}_{\frac{1}{2},k}\bigl(\mathbb{B}(\mathcal{T}),\Upsilon\bigr)\bigr|=
\sum_{X\in\mathcal{C}_{\frac{1}{2},k}(\mathbb{B}(\mathcal{T}),\Upsilon):\ X>\hat{0}}
\mu_{\mathcal{C}}(\hat{0},X)\\
\cdot\sum_{z\in\mathcal{E}(\mathbb{B}(\mathcal{T}),X):\ z>\hat{0}}\mu_{\mathcal{E}}(\hat{0},z)\cdot\binom{|\mathcal{T}|-\rho(z)}{|\mathcal{T}|-k}\ .
\end{multline}
Restate~(\ref{eq:27}) in the following way:

\begin{proposition}
\label{prop:3}
The number $\#\mathbf{K}^{\ast}_k(\mathcal{M})$ of tope committees, of cardinality~$k$, $3\leq k\leq
|\mathcal{T}|-3$, for the oriented matroid $\mathcal{M}:=(E_t,\mathcal{T})$, is
\begin{multline*}
\#\mathbf{K}^{\ast}_k(\mathcal{M})=\#\mathbf{K}^{\ast}_{|\mathcal{T}|-k}(\mathcal{M})=
\sum_{
\mathcal{G}\in\{\ \{\bigcup_{e\in E}\binom{\mathcal{T}^+_e}{\lceil(\ell+1)/2\rceil}\}:\ E\subseteq E_t,\ |E|>0\ \}
}\mu_{\boldsymbol{\mathcal{C}}}(\hat{0},\mathcal{G})\\
\cdot\sum_{G\in\boldsymbol{\mathcal{E}}(\mathcal{G}):\ 0<|G|\leq \ell}\mu_{\boldsymbol{\mathcal{E}}}(\hat{0},G)\cdot\binom{|\mathcal{T}|-|G|}{\ell-|G|}\ ,
\end{multline*}
where $\ell\in\{k,|\mathcal{T}|-k\}$;
$\mu_{\boldsymbol{\mathcal{C}}}(\cdot,\cdot)$ denotes the M\"{o}bius function of the family~$\boldsymbol{\mathcal{C}}$ $:=\{\hat{0}\}\dot\cup\bigl\{\{\bigcup_{e\in E}\binom{\mathcal{T}^+_e}{\lceil(\ell+1)/2\rceil}\}:\ E\subseteq E_t,\ |E|>0\bigr\}$ ordered by inclusion.
\end{proposition}

\section{Convex Sets and Tope Committees}

\label{s:2}

Let the antichain $\Upsilon:=\{\upsilon_1,\ldots,\upsilon_t\}\subset\mathbb{B}(\mathcal{T})^{(|\mathcal{T}|/2)}$  again represent the family of positive halspaces of the oriented matroid $\mathcal{M}$ in the Boolean lattice~$\mathbb{B}(\mathcal{T})$ of all subsets of the tope set $\mathcal{T}$. We have
\begin{multline}
\label{eq:1}
\bigl|\mathbf{I}_{\frac{1}{2},k}\bigl(\mathbb{B}(\mathcal{T}),\Upsilon\bigr)\bigr|
=\sum_{D\subseteq \mathbb{B}(\mathcal{T})^{(\lceil(k+1)/2\rceil)}\cap\mathfrak{I}(\Upsilon):\ |D|>0}
(-1)^{|D|}\;\\ \cdot\left(\sum_{C\subseteq\Upsilon:\ D\subseteq\mathfrak{I}(C)}(-1)^{|C|}\right)\!\!\!
\binom{|\mathcal{T}|-\rho(\vee_{d\in D}d)}{|\mathcal{T}|-k}\ ,\ \ \ 3\leq k\leq|\mathcal{T}|-3\ ,
\end{multline}
cf.~(\ref{eq:17}).

Consider the mapping
\begin{equation}
\label{eq:3}
\begin{split}
\gamma_k:\ \mathbb{B}(\mathcal{T})^{(\lceil(k+1)/2\rceil)}\cap\mathfrak{I}(\Upsilon)
&\to L_{\conv}(\mathcal{M})\ ,\\  d&\mapsto\bmax\{A\in L_{\conv}(\mathcal{M}):\
d\subseteq\mathcal{T}^+_A\}\ ,
\end{split}
\end{equation}
that sends a $(\lceil(k+1)/2\rceil)$-subset of topes $d\in\mathfrak{I}(\Upsilon)$ to the inclusion-maximum convex subset $A\subset E_t$ with the property $d\subseteq\mathcal{T}^+_A$;
we are actually interested in such a mapping to the subposet $L_{\conv,\ \geq \lceil(k+1)/2\rceil}(\mathcal{M})$,
the order ideal of the semilattice $L_{\conv}(\mathcal{M})$ defined as
$L_{\conv,\ \geq \lceil(k+1)/2\rceil}(\mathcal{M}):=\bigl\{A$ $\in L_{\conv}(\mathcal{M}):\ |\mathcal{T}^+_A|\geq\lceil(k+1)/2\rceil\bigr\}$.

Fix a nonempty subset $D\subseteq \mathbb{B}(\mathcal{T})^{(\lceil(k+1)/2\rceil)}\cap\mathfrak{I}(\Upsilon)$ and consider the
{\em blocker\/}~$\mathcal{B}\bigl(\gamma_k(D)\bigr)$ of the image $\gamma_k(D)$; if we let $\bminit\gamma_k(D)$ denote
the subfamily of all inclusion-minimal sets from the family $\gamma_k(D)$ then $\mathcal{B}\bigl(\gamma_k(D)\bigr)$ $=
\mathcal{B}\bigl(\bminit \gamma_k(D)\bigr)$.

Let $\Delta^{\ast}(D)$ be the abstract simplicial complex whose facets are the complements $E_t-B$ of the sets $B\in\mathcal{B}\bigl(\bminit\gamma_k(D)\bigr)$ from the blocker of the {\em Sperner family\/} $\bminit\gamma_k(D)$, and let
$\Delta(D)$ be the complex
whose facets are the complements $E_t-G$ of the sets $G\in\bminit\gamma_k(D)$; if the complexes~$\Delta(D)$ and~$\Delta^{\ast}(D)$ have the same vertex set then $\Delta^{\ast}(D)$ is the {\em Alexander dual\/} of~$\Delta(D)$. The {\em reduced Euler characteristics\/} $\widetilde{\chi}(\cdot)$ of the complexes satisfy the equality $\widetilde{\chi}\bigl(\Delta^{\ast}(D)\bigr)=(-1)^{t-1}
\widetilde{\chi}\bigl(\Delta(D)\bigr)$.

For a subset $C:=(\upsilon_{i_1},\ldots,\upsilon_{i_j})\subseteq\Upsilon$ we have $D\subseteq\mathfrak{I}(C)$ iff the collection of indices $\{i_1,\ldots,i_j\}$ is a blocking set for the family $\bminit\gamma_k(D)$; therefore
\begin{equation*}
\sum_{C\subseteq\Upsilon:\ D\subseteq\mathfrak{I}(C)}(-1)^{|C|}=(-1)^{t-1}\widetilde{\chi}\bigl(\Delta^{\ast}(D)\bigr)\ .
\end{equation*}
If $\bigcup_{F\in\bminit\gamma_k(D)}F\neq E_t$ then the complex $\Delta^{\ast}(D)$ is a {\em cone\/} and, as a consequence,
$\widetilde{\chi}\bigl(\Delta^{\ast}(D)\bigr)=0$.

Rewrite~(\ref{eq:1}) in the following way:
\begin{multline}
\label{eq:2}
\bigl|\mathbf{I}_{\frac{1}{2},k}\bigl(\mathbb{B}(\mathcal{T}),\Upsilon\bigr)\bigr|=
\sum_{\substack{D\subseteq \mathbb{B}(\mathcal{T})^{(\lceil(k+1)/2\rceil)}\cap\mathfrak{I}(\Upsilon):\\
\bigcup_{F\in\bminit\gamma_k(D)}F=E_t}}
(-1)^{|D|}\\ \cdot\widetilde{\chi}\bigl(\Delta(D)\bigr)\cdot
\binom{|\mathcal{T}|-\rho(\vee_{d\in D}d)}{k-\rho(\vee_{d\in D}d)}\ ;
\end{multline}
note that singleton sets $D:=\{d\}$, where $d\in\mathbb{B}(\mathcal{T})^{(\lceil(k+1)/2\rceil)}\cap\mathfrak{I}(\Upsilon)$, do not play a r\^{o}le in~(\ref{eq:2}).

Given a subset $D\subseteq \mathbb{B}(\mathcal{T})^{(\lceil(k+1)/2\rceil)}\cap\mathfrak{I}(\Upsilon)$ such that $\bigcup_{F\in\bminit\gamma_k(D)}F= E_t$, let\hfill $\mathcal{S}(D)$\hfill denote\hfill the\hfill family\hfill of\hfill the\hfill unions\hfill $\{\bigcup_{F\in\mathcal{F}}F:\ \mathcal{F}\subseteq\bminit\gamma_k(D),$\\ $\#\mathcal{F}>0\}$ ordered by inclusion and augmented by a new least element $\hat{0}$; the greatest element $\hat{1}$ of the lattice
$\mathcal{S}(D)$ is the ground set $E_t$. The reduced Euler characteristic $\widetilde{\chi}\bigl(\Delta(D)\bigr)=\sum_{\substack{\mathcal{F}\subseteq\bminit\gamma_k(D):\\ \bigcup_{F\in\mathcal{F}}F=E_t}}(-1)^{\#\mathcal{F}}$ of the complex $\Delta(D)$ is equal
to the {\em M\"{o}bius number\/}~$\mu_{\mathcal{S}(D)}(\hat{0},\hat{1})$ and, in particular, to~$(-1)^{\#\bminit\gamma_k(D)}$ when the sets in the family $\bminit\gamma_k(D)$ are pairwise disjoint.
Restate observation~(\ref{eq:2}):
\begin{proposition}
\label{prop:1}
The number $\#\mathbf{K}^{\ast}_k(\mathcal{M})$ of tope committees, of cardinality~$k$, $3\leq k\leq
|\mathcal{T}|-3$, for the oriented matroid $\mathcal{M}:=(E_t,\mathcal{T})$, is
\begin{multline}
\label{eq:7}
\#\mathbf{K}^{\ast}_k(\mathcal{M})=\#\mathbf{K}^{\ast}_{|\mathcal{T}|-k}(\mathcal{M})\\=
\sum_{\substack{\mathcal{G}\subseteq\bigcup_{e\in E_t}\binom{\mathcal{T}^+_e}{\lceil(\ell+1)/2\rceil}:\\
1<\#\mathcal{G}\leq\binom{\ell}{\lceil(\ell+1)/2\rceil},\\
\bigcup_{F\in\bminit\gamma_{\ell}(\mathcal{G})}F=E_t,\
|\bigcup_{G\in\mathcal{G}}G|\leq\ell
}}
(-1)^{\#\mathcal{G}}\cdot\mu_{\mathcal{S}(\mathcal{G})}(\hat{0},\hat{1})\cdot
\binom{|\mathcal{T}|-|\bigcup_{G\in\mathcal{G}}G|}{\ell-|\bigcup_{G\in\mathcal{G}}G|}\ ,
\end{multline}
where $\ell\in\{k,|\mathcal{T}|-k\}$.
\end{proposition}

Consider the abstract simplicial complex whose facets are the positive halfspaces of the oriented matroid $\mathcal{M}$. If some its relevant $(\lceil(k+1)/2\rceil-1)$-dimensional faces, sets from the family
$\bigcup_{e\in E_t}\binom{\mathcal{T}^+_e}{\lceil(k+1)/2\rceil}$, are {\em free\/} --- each of them is  contained in exactly one facet $\mathcal{T}^+_e$, for some element $e\in E_t$ --- then the M\"{o}bius numbers~$\mu_{\mathcal{S}(\mathcal{G})}(\hat{0},\hat{1})$ in~(\ref{eq:7}), under $\ell:=k$, are all equal to~$(-1)^t$:

\begin{corollary}
Let $k$ be an integer, $3\leq k\leq
|\mathcal{T}|-3$. If for any family $\mathcal{G}\subseteq\bigcup_{e\in E_t}\binom{\mathcal{T}^+_e}{\lceil(k+1)/2\rceil}$ such that $\bigcup_{F\in\bminit\gamma_k(\mathcal{G})}F=E_t$ and $|\bigcup_{G\in\mathcal{G}}G|\leq k$,
it holds $|\gamma_k(G)|=1$, for any set $G\in\mathcal{G}$, then
the number $\#\mathbf{K}^{\ast}_k(\mathcal{M})$ of tope committees, of cardinality $k$,  for the oriented matroid $\mathcal{M}:=(E_t,\mathcal{T})$, is
\begin{multline}
\label{eq:11}
\#\mathbf{K}^{\ast}_k(\mathcal{M})=\#\mathbf{K}^{\ast}_{|\mathcal{T}|-k}(\mathcal{M})\\=(-1)^t\cdot
\sum_{\substack{\mathcal{G}\subseteq\bigcup_{e\in E_t}\binom{\mathcal{T}^+_e}{\lceil(k+1)/2\rceil}:\\
t\leq\#\mathcal{G}\leq\binom{k}{\lceil(k+1)/2\rceil},\\
\bigcup_{F\in\bminit\gamma_k(\mathcal{G})}F=E_t,\ |\bigcup_{G\in\mathcal{G}}G|\leq k
}}
(-1)^{\#\mathcal{G}}\cdot
\binom{|\mathcal{T}|-|\bigcup_{G\in\mathcal{G}}G|}{k-|\bigcup_{G\in\mathcal{G}}G|}\ .
\end{multline}
\end{corollary}

\section{Relative Blocking in~Posets~Isomorphic to the~Face~Lattices
of~Crosspolytopes}

\label{s:5}

Consider a poset $\mathbf{O}'(m)$, with the rank function $\rho(\cdot)$, which is isomorphic to the graded face meet-semilattice of the boundary of a $m$-dimensional {\em crosspolytope} and is defined in the following way: the semilattice~$\mathbf{O}'(m)$ is composed of all subsets, free of opposites, of a set~$\{-m,$ $\ldots,-1,1,\ldots,m\}$, ordered by inclusion. We denote by $\mathbf{O}(m)$ the lattice $\mathbf{O}(m):=\mathbf{O}'(m)$ $\dot\cup\{\hat{1}\}$, where~$\hat{1}$ is a new greatest element.
Let~$\varLambda\subset\mathbf{O}'(m)$ be a nontrivial antichain in the lattice $\mathbf{O}(m)$.

For a rational number $r$, $0\leq r<1$, and for a positive integer number~$k$, we define the set $\mathbf{I}_{r,k}\bigl(\mathbf{O}'(m),
\varLambda\bigr)$ of {\em relatively $r$-blocking elements}, of rank~$k$, for the antichain~$\varLambda$ in analogy with the sets $\mathbf{I}_{r,k}\bigl(\mathbb{B}(n),\cdot\bigr)$ for antichains in Boolean lattices, cf.~(\ref{eq:12}):
\begin{multline*}
\mathbf{I}_{r,k}\bigl(\mathbf{O}'(m),
\varLambda\bigr):=\bigl\{b\in\mathbf{O}'(m):\ \rho(b)=k,\\
\rho(b\wedge\lambda)>r\cdot k\ \ \ \forall
\lambda\in\varLambda\bigr\}\subset\mathbf{O}'(m)^{(k)}\ ,
\end{multline*}
where $\mathbf{O}'(m)^{(k)}$ is the $k$th layer of the semilattice $\mathbf{O}'(m)$.

On the one hand, we have
\begin{multline}
\label{eq:20}
\bigl|\mathbf{I}_{r,k}\bigl(\mathbf{O}'(m),\varLambda\bigr)\bigr|=
\sum_{X\in\overset{\circ}{\mathcal{C}}_{r,k}(\mathbf{O}'(m),\varLambda):\ X>\hat{0}}
\mu_{\overset{\circ}{\mathcal{C}}}(\hat{0},X)\\
\cdot\sum_{z\in\overset{\circ}{\mathcal{E}}(\mathbf{O}'(m),X):\ z>\hat{0}}\mu_{\overset{\circ}{\mathcal{E}}}(\hat{0},z)\cdot 2^{k-\rho(z)}\cdot\binom{m-\rho(z)}{m-k}\ ,
\end{multline}
cf.~(\ref{eq:19}), where $\overset{\circ}{\mathcal{C}}_{r,k}(\mathbf{O}'(m),\varLambda)$ denotes the join-semilattice of all sets
from the family $\{\mathbf{O}'(m)^{(\lfloor r\cdot k\rfloor+1)}\cap\mathfrak{I}(C):\ C\subseteq \varLambda,\ |C|>0\}$
ordered by inclusion and augmented by a new least element $\hat{0}$; the greatest element~$\hat{1}$ of the lattice $\overset{\circ}{\mathcal{C}}_{r,k}(\mathbf{O}'(m),\varLambda)$ is the set $\mathbf{O}'(m)^{(\lfloor r\cdot k\rfloor+1)}\cap\mathfrak{I}(\varLambda)$. For an element $X\in\overset{\circ}{\mathcal{C}}_{r,k}(\mathbf{O}'(m),\varLambda)$, the notation $\overset{\circ}{\mathcal{E}}(\mathbf{O}'(m),X)$ is used to denote the sub-join-semilattice of the lattice $\mathbf{O}(m)$ generated by the set $X\subset\mathbf{O}'(m)$, with the greatest element of~$\mathbf{O}(m)$ deleted from it, and augmented by a new least element $\hat{0}$. The M\"{o}bius functions of the posets~$\overset{\circ}{\mathcal{C}}_{r,k}(\mathbf{O}'(m),\varLambda)$ and~$\overset{\circ}{\mathcal{E}}(\mathbf{O}'(m),X)$ are denoted by~$\mu_{\overset{\circ}{\mathcal{C}}}(\cdot,\cdot)$ and~$\mu_{\overset{\circ}{\mathcal{E}}}(\cdot,\cdot)$, respectively;  $\rho(\cdot)$ denotes the rank of an element in the poset~$\mathbf{O}'(m)$.

On the other hand, we have
\begin{multline}
\label{eq:18}
\bigl|\mathbf{I}_{r,k}\bigl(\mathbf{O}'(m),\varLambda\bigr)\bigr|
=\sum_{\substack{D\subseteq \mathbf{O}'(m)^{(\lfloor r\cdot k\rfloor+1)}\cap\mathfrak{I}(\varLambda):\\ |D|>0,\ \bigvee_{d\in D}d\neq\hat{1}}}(-1)^{|D|}\;\\ \cdot\left(\sum_{C\subseteq\varLambda:\ D\subseteq\mathfrak{I}(C)}(-1)^{|C|}\right)\cdot 2^{k-\rho(\vee_{d\in D}d)}\cdot
\binom{m-\rho(\vee_{d\in D}d)}{m-k}\ ,
\end{multline}
cf.~(\ref{eq:17}).

\section{Tope Committees Containing no Pairs of Opposites}

\label{s:3}

Counting tope committees, that are free of opposites, for the oriented matroid $\mathcal{M}$, we follow the reasoning scheme from Section~\ref{s:2}, but we work now with the family~$\mathbf{O}'(\mathcal{T})$ of tope subsets that are free of opposites and ordered by inclusion; the semilattice~$\mathbf{O}'(\mathcal{T})$ is isomorphic to the face poset of the boundary of a crosspolytope of dimension $|\mathcal{T}|/2$. The lattice $\mathbf{O}(\mathcal{T})$ $:=\mathbf{O}'(\mathcal{T})\dot\cup\{\hat{1}\}$ is the semilattice~$\mathbf{O}'(\mathcal{T})$ augmented by a new greatest element $\hat{1}$. We again turn
to the mapping~$\gamma_k:\ \mathbb{B}(\mathcal{T})^{(\lceil(k+1)/2\rceil)}\cap\mathfrak{I}(\Upsilon)$ $=\mathbf{O}'(\mathcal{T})^{(\lceil(k+1)/2\rceil)}\cap\mathfrak{I}(\Upsilon)\to L_{\conv}(\mathcal{M})$ defined in~(\ref{eq:3}), and to the lattices~$\mathcal{S}(\cdot)$ considered in Section~\ref{s:2}.

If $\mathcal{G}$ is a family of tope subsets which are free of opposites then we denote by~$\overset{\circ}{\boldsymbol{\mathcal{E}}}(\mathcal{G})$ the join-semilattice $\{\bigcup_{F\in\mathcal{F}}F:\ \mathcal{F}\subseteq\mathcal{G},\ \#\mathcal{F}>0,\ \text{$\bigcup_{F\in\mathcal{F}}F$ free}$ $\text{of opposites}\}$ composed of the unions, free of opposites, of the sets from the family $\mathcal{G}$ ordered by inclusion and augmented by a new least element $\hat{0}$;
the M\"{o}bius function of the poset $\overset{\circ}{\boldsymbol{\mathcal{E}}}(\mathcal{G})$ is denoted by $\mu_{\overset{\circ}{\boldsymbol{\mathcal{E}}}}(\cdot,\cdot)$.

Formula~(\ref{eq:21}) below is deduced from~(\ref{eq:20}). Formulas~(\ref{eq:9}) and~(\ref{eq:10}) are deduced from~(\ref{eq:18}); they are direct analogues of formulas~(\ref{eq:7}) and~(\ref{eq:11}), respectively. See also~\cite[Section~3]{M-Reorientations}.

\begin{theorem}
The number $\#\overset{\circ}{\mathbf{K}}{}^{\ast}_k(\mathcal{M})$ of tope committees which are free of opposites, of cardinality $k$, $3\leq k\leq|\mathcal{T}|/2$, for the oriented matroid~$\mathcal{M}$ $:=(E_t,\mathcal{T})$, is:

{\rm(i)}

\nopagebreak[4]

\begin{multline}
\label{eq:21}
\#\overset{\circ}{\mathbf{K}}{}^{\ast}_k(\mathcal{M})=
\sum_{
\mathcal{G}\in\{\ \{\bigcup_{e\in E}\binom{\mathcal{T}^+_e}{\lceil(k+1)/2\rceil}\}:\ E\subseteq E_t,\ |E|>0\ \}
}\mu_{\overset{\circ}{\boldsymbol{\mathcal{C}}}}(\hat{0},\mathcal{G})\\
\cdot\sum_{G\in\overset{\circ}{\boldsymbol{\mathcal{E}}}(\mathcal{G}):\ 0<|G|\leq k}\mu_{\overset{\circ}{\boldsymbol{\mathcal{E}}}}(\hat{0},G)\cdot
2^{{k-|G|}}\cdot\binom{\frac{1}{2}|\mathcal{T}|-|G|}{k-|G|}\ ,
\end{multline}
where\hfill
$\mu_{\overset{\circ}{\boldsymbol{\mathcal{C}}}}(\cdot,\cdot)$\hfill denotes\hfill the\hfill M\"{o}bius\hfill function\hfill of\hfill the\hfill family\hfill  $\overset{\circ}{\boldsymbol{\mathcal{C}}}:=\{\hat{0}\}$ \newline $\dot\cup\bigl\{\{\bigcup_{e\in E}\binom{\mathcal{T}^+_e}{\lceil(k+1)/2\rceil}\}:\ E\subseteq E_t,\ |E|>0\bigr\}$ ordered by inclusion.

{\rm(ii)}

\begin{multline}
\label{eq:9}
\#\overset{\circ}{\mathbf{K}}{}^{\ast}_k(\mathcal{M})=
\sum_{\substack{\mathcal{G}\subseteq\bigcup_{e\in E_t}\binom{\mathcal{T}^+_e}{\lceil(k+1)/2\rceil}:\\
1<\#\mathcal{G}\leq\binom{k}{\lceil(k+1)/2\rceil},\\
\bigcup_{G\in\mathcal{G}}G\text{\rm\ free of opposites},\\
\bigcup_{F\in\bminit\gamma_k(\mathcal{G})}F=E_t,\ |\bigcup_{G\in\mathcal{G}}G|\leq k
}}
(-1)^{\#\mathcal{G}}\cdot\mu_{\mathcal{S}(\mathcal{G})}(\hat{0},\hat{1})\\ \cdot
2^{{k-|\bigcup_{G\in\mathcal{G}}G|}}\cdot
\binom{\frac{1}{2}|\mathcal{T}|-|\bigcup_{G\in\mathcal{G}}G|}{k-|\bigcup_{G\in\mathcal{G}}G|}\ .
\end{multline}
\noindent In particular, if for any family $\mathcal{G}\subseteq\bigcup_{e\in E_t}\binom{\mathcal{T}^+_e}{\lceil(k+1)/2\rceil}$ such that $\bigcup_{G\in\mathcal{G}}G$ is free of opposites, $\bigcup_{F\in\bminit\gamma_k(\mathcal{G})}F$ $=E_t$ and $|\bigcup_{G\in\mathcal{G}}G|\leq k$,
it holds $|\gamma_k(G)|$ $=1$, for any set $G\in\mathcal{G}$, then
\begin{multline}
\label{eq:10}
\overset{\circ}{\mathbf{K}}{}^{\ast}_k(\mathcal{M})=(-1)^t\\ \cdot
\sum_{\substack{\mathcal{G}\subseteq\bigcup_{e\in E_t}\binom{\mathcal{T}^+_e}{\lceil(k+1)/2\rceil}:\\
t\leq\#\mathcal{G}\leq\binom{k}{\lceil(k+1)/2\rceil},\\
\bigcup_{G\in\mathcal{G}}G\text{\rm\ free of opposites},\\
\bigcup_{F\in\bminit\gamma_k(\mathcal{G})}F=E_t,\ |\bigcup_{G\in\mathcal{G}}G|\leq k}}
(-1)^{\#\mathcal{G}}\cdot
2^{{k-|\bigcup_{G\in\mathcal{G}}G|}}\cdot
\binom{\frac{1}{2}|\mathcal{T}|-|\bigcup_{G\in\mathcal{G}}G|}{k-|\bigcup_{G\in\mathcal{G}}G|}\ .
\end{multline}
\end{theorem}

\section{Relative Blocking in Principal Order Ideals of Binomial Posets}

\label{s:6}

In this section we mention an analogue of formula~(\ref{eq:2}) in the more general context of binomial posets.

Let $P$ be a graded lattice of rank $n$ which is a principal order ideal of some {\em binomial poset}. The {\em factorial function\/} $\mathrm{B}(k)$ of $P$ counts the number of {\em maximal chains\/} in any {\em interval\/} of length $k$ in $P$. The number $\left[\begin{smallmatrix}j\\ i\end{smallmatrix}\right]$ of elements of rank $i$ in an interval of length $j$ in $P$ is equal to~$\tfrac{\mathrm{B}(j)}{\mathrm{B}(i)\cdot\mathrm{B}(j-i)}$, see~\cite[\S{}3.15]{St1}.

Let~$\varLambda$ be a nontrivial antichain in the lattice $P$. If $r$ is a rational number, $0\leq r<1$, and $k$ is a positive integer, then the set $\mathbf{I}_{r,k}(P,\varLambda)$ of {\em relatively $r$-blocking elements}, of rank $k$, for the antichain $\varLambda$ in $P$, is defined as follows:
\begin{equation*}
\mathbf{I}_{r,k}\bigl(P,
\varLambda\bigr):=\bigl\{b\in P:\ \rho(b)=k,\
\rho(b\wedge\lambda)>r\cdot k\ \ \ \forall
\lambda\in\varLambda\bigr\}\subset P^{(k)}\ ,
\end{equation*}
where $\rho(\cdot)$ is the rank function of $P$, and $P^{(k)}$ is the $k$th layer of $P$.

Let $\mathcal{N}(\varLambda)$ be the abstract simplicial complex whose facets are the inclusion-maximal sets of indices $\{i_1,\ldots,i_j\}$ such that for the corresponding antichains $\{\lambda_{i_1},\ldots,\lambda_{i_j}\}\subseteq\varLambda$ it holds $\lambda_{i_1}\wedge\cdots\wedge\lambda_{i_j}>\hat{0}$, where $\hat{0}$ is the least element of $P$. If the poset $P$ is the Boolean lattice $\mathbb{B}(n)$ then the complex~$\mathcal{N}(\varLambda)$ is the {\em nerve\/} of the corresponding Sperner family; see, e.g.,~\cite[\S{}10]{B-TM} on the topological combinatorics of the nerve.

Set $\nu(r\cdot k):=\lfloor r\cdot k\rfloor+1$. Let
\begin{align*}
\mathfrak{c}_{r,k}:P^{(\nu(r\cdot k))}\cap\mathfrak{I}(\varLambda)&\to\mathcal{N}(\varLambda)\ ,\\
d&\mapsto\bmax\bigl\{N\in\mathcal{N}(\varLambda):\ d\leq\bigwedge_{i\in N}\lambda_i\bigr\}
\end{align*}
be the mapping that reflects an element $d$, of rank $\nu(r\cdot k)$, of the order ideal~$\mathfrak{I}(\varLambda)$ generated by the antichain $\varLambda$ to the inclusion-maximum face of the complex~$\mathcal{N}(\varLambda)$ with the property $d\leq\bigwedge_{i\in N}\lambda_i$.

Associate\hfill to\hfill a\hfill subset\hfill $D\subseteq P^{(\nu(r\cdot k))}\cap\mathfrak{I}(\varLambda)$,\hfill such\hfill that\hfill $|\bigcup_{F\in\bminit\mathfrak{c}_{r,k}(D)}F|$\\ $=|\varLambda|$,\hfill a\hfill poset\hfill $\mathcal{S}(D)$\hfill which\hfill is\hfill the\hfill family\hfill $\{\bigcup_{F\in\mathcal{F}}F:\ \mathcal{F}\subseteq\bminit\mathfrak{c}_{r,k}(D),\#\mathcal{F}$\\ $>0\}$\hfill ordered\hfill by\hfill inclusion,\hfill with\hfill a\hfill new\hfill least\hfill element\hfill $\hat{0}$\hfill adjoined;\\ here~$\bminit\mathfrak{c}_{r,k}(D)$ denotes the subfamily of all inclusion-minimal sets from the image $\mathfrak{c}_{r,k}(D)$. Let $\mu_{\mathcal{S}(D)}(\hat{0},\hat{1})$ denote the corresponding M\"{o}bius number, where $\hat{1}$ is the greatest element of $\mathcal{S}(D)$.

Suppose that $\nu(r\cdot k)\leq \min_{\lambda\in\varLambda}\rho(\lambda)$.
Since
\begin{multline*}
\bigl|\mathbf{I}_{r,k}\bigl(P,\varLambda\bigr)\bigr|
=\sum_{D\subseteq P^{(\nu(r\cdot k))}\cap\mathfrak{I}(\varLambda):\ |D|>0}
(-1)^{|D|}\;\\ \cdot\left(\sum_{C\subseteq\varLambda:\ D\subseteq\mathfrak{I}(C)}(-1)^{|C|}\right)\!\!\!
\begin{bmatrix}n-\rho(\vee_{d\in D}d)\\n-k\end{bmatrix}\ ,
\end{multline*}
by~\cite[(5.4)]{M-Relative}, we have
\begin{multline*}
\bigl|\mathbf{I}_{r,k}\bigl(P,\varLambda\bigr)\bigr|
=\sum_{\substack{D\subseteq P^{(\nu(r\cdot k))}\cap\mathfrak{I}(\varLambda):\\
1\leq|D|\leq\left[\begin{smallmatrix}k\\ \nu(r\cdot k)\end{smallmatrix}\right],\\
|\bigcup_{F\in\bminit\mathfrak{c}_{r,k}(D)}F|=|\varLambda|,\\
\rho(\vee_{d\in D}d)\leq k}}
(-1)^{|D|}\cdot\mu_{\mathcal{S}(D)}(\hat{0},\hat{1})\cdot
\begin{bmatrix}n-\rho(\vee_{d\in D}d)\\k-\rho(\vee_{d\in D}d)\end{bmatrix}\ .
\end{multline*}


\end{document}